\documentclass[a4paper,reqno]{amsart}

\usepackage{amsmath}
\usepackage{amssymb}
\usepackage{amscd}
\usepackage{amsthm}
\usepackage{type1cm}
\usepackage{tcolorbox}
\usepackage{mathrsfs}

\usepackage[colorlinks,citecolor=darkgreen,linkcolor=red]{hyperref}
\definecolor{darkgreen}{rgb}{0.0, 0.6, 0.0}

\usepackage[all]{xy}
\input xy
\xyoption{all}
\tcbuselibrary{breakable, skins, theorems}

\usepackage{tikz-cd}

\def\C{\mathcal{C}}
\def\D{\mathcal{D}}
\def\E{\mathcal{E}}
\def\H{\mathcal{H}}
\def\K{\mathcal{K}}

\def\O{\mathcal{O}}
\def\T{\mathcal{T}}
\def\U{\mathcal{U}}

\DeclareMathOperator{\Mod}{\mathsf{Mod}}
\DeclareMathOperator{\md}{\mathsf{mod}}
\renewcommand{\mod}{\md}
\DeclareMathOperator{\fl}{\mathsf{fl}}
\DeclareMathOperator{\proj}{\mathsf{proj}}
\DeclareMathOperator{\refl}{\mathsf{ref}}
\DeclareMathOperator{\per}{\mathsf{per}}
\DeclareMathOperator{\pvd}{\mathsf{pvd}}

\DeclareMathOperator{\add}{\mathsf{add}}
\DeclareMathOperator{\Add}{\mathsf{Add}}
\DeclareMathOperator{\thick}{\mathsf{thick}}
\DeclareMathOperator{\qmod}{\mathsf{qmod}}
\DeclareMathOperator{\Coh}{\mathsf{Coh}}

\DeclareMathOperator{\End}{End}
\DeclareMathOperator{\Ext}{Ext}
\DeclareMathOperator{\op}{op}

\DeclareMathOperator{\Ker}{Ker}
\DeclareMathOperator{\Cok}{Cok}

\DeclareMathOperator{\rad}{rad}
\DeclareMathOperator{\tp}{top}
\renewcommand{\top}{\tp}
\DeclareMathOperator{\RHom}{\mathbb{R}Hom}
\DeclareMathOperator{\REnd}{\mathbb{R}End}
\DeclareMathOperator{\silt}{silt}
\DeclareMathOperator{\ctilt}{ctilt}
\DeclareMathOperator{\holim}{holim}

\def\gl{\mathop{\rm gl.dim}\nolimits}

\def\rgl{\mathop{\rm r.gl.dim}\nolimits}
\def\wgl{\mathop{\rm w.gl.dim}\nolimits}
\def\pd{\mathop{\rm proj.dim}\nolimits}
\def\id{\mathop{\rm inj.dim}\nolimits}
\def\fd{\mathop{\rm flat.dim}\nolimits}

\theoremstyle{definition}
\newtheorem{Thm}{Theorem}[section]
\newtheorem{Lem}[Thm]{Lemma}
\newtheorem{Prop}[Thm]{Proposition}
\newtheorem{Cor}[Thm]{Corollary}
\newtheorem{Def}[Thm]{Definition}
\newtheorem{Ex}[Thm]{Example}
\newtheorem{Rem}[Thm]{Remark}
\newtheorem{Conj}[Thm]{Conjecture}
\newtheorem{Propdef}[Thm]{Proposition-Definition}
\newtheorem{Ques}[Thm]{Question}

\newcommand{\FRAC}[2]{\leavevmode\kern.1em\raise.5ex\hbox{\the\scriptfont0 #1}\kern-.1em/\kern-.15em\lower.25ex\hbox{\the\scriptfont0 #2}}

\title{On silting mutations preserving global dimension}
\author{Ryu Tomonaga}
\address{Graduate School of Mathematical Sciences, The University of Tokyo, 3-8-1 Komaba, Meguro-ku, Tokyo, 153-8914, Japan}
\email{ryu-tomonaga@g.ecc.u-tokyo.ac.jp}

\begin{document}
\begin{abstract}
A $d$-silting object is a silting object whose derived endomorphism algebra has global dimension $d$ or less. We give an equivalent condition, which can be stated in terms of dg quivers, for silting mutations to preserve the $d$-siltingness under a mild assumption. Moreover, we show that this mild assumption is always satisfied by $\nu_d$-finite algebras.

As an application, we give counterexamples to the open question by Herschend-Iyama-Oppermann: the quivers of higher hereditary algebras are acyclic. Our examples consist of a $2$-representation tame algebra with a $2$-cycle, and a $3$-homogeneous $2$-representation finite algebra with a cycle.
\end{abstract}

\maketitle
\tableofcontents

\section*{Introduction}

The notion of tilting objects is indispensable to construct derived equivalences and gives a deep connections among many areas of mathematics such as representation theory, algebraic geometry and mathematical physics. The notion of silting objects is a natural generalization of that of tilting objects from the view point of mutation \cite{AI}, which is parallel to that the notion of connective dg algebras is a generalization of that of algebras. Silting mutation is a fundamental way to reproduce silting objects from a given one and has a strong relationship with mutations in cluster algebras \cite{AIR14}.

On the other hand, the notion of global dimension is a fundamental invariant of algebras which measures how complex the module category or the derived category is and plays an essential role in higher Auslander-Reiten theory. To deal with global dimension systematically, in \cite{HI}, the notion of $d$-silting objects is introduced for $d\geq1$: a $d$-silting object is a silting object whose derived endomorphism algebra has global dimension $d$ or less. This is a natural generalization of the notion of $d$-tilting objects which is extensively studied in \cite{BHI,HIMO,Tom25b,Tom25a}. In \cite{HI}, they establish connections called silting correspondence between $d$-silting objects and cluster tilting objects or silting objects of its $(d+1)$-Calabi-Yau completion \cite{Kel11}.

In this paper, we investigate when the silting mutation preserves $d$-siltingness. The following theorem gives a clear answer.

\begin{Thm}[Corollary \ref{easychar}]
Let $A$ be a proper connective dg algebra with $\gl A\leq d$. Write $A=P\oplus Q$ with $(\add P)\cap(\add Q)=0$, and set $S:=\top H^0P$. Then the mutation $\mu_P^-(A)$ is $d$-silting if and only if
\[
 \D(A)(S,Q[d])=0.
\]
Moreover, if $\nu$ denotes the Serre functor of $\per A$ and $\nu_d:=\nu[-d]$, then the following conditions are equivalent:
\begin{enumerate}
\item $\mu_P^-(A)\geq\nu_d^{-1}A$;
\item $\pd_AS<d$;
\item $\mu_P^-(A)$ is $d$-silting and $\D(A)(S,S[d])=0$.
\end{enumerate}
\end{Thm}

The criterion becomes combinatorial for dg path algebras. Let $A=kQ$ be proper, where $Q$ is a finite graded quiver with $Q_1^{>0}=Q_1^{\leq-d}=\emptyset$, and assume that the differential of every arrow lies in $kQ_{\geq2}$. Note that every proper connective dg algebra over an algebraically closed field with global dimension at most $d$ has such a description. We compute extensions between simple dg modules, and more generally the homological dimensions of $A$, from $Q$ (Theorem \ref{extsimpath} and Corollary \ref{gldimpath}). In particular, we obtain the following criterion.

\begin{Thm}[Dg-quiver criterion; Corollary~\ref{easycharpath}]
Let $i\in Q_0$ and assume that there is no loop of degree $-d+1$ at $i$. Then $\mu_{e_iA}^-(A)$ is $d$-silting if and only if no arrow of degree $-d+1$ ends at $i$.
\end{Thm}

For dg path algebras, the same criterion can be recovered from the explicit mutation procedure of \cite{Opp}. Our argument instead explains the criterion intrinsically and applies without choosing a dg-quiver presentation.

The self-extension hypothesis in the preceding theorem is automatic under $\nu_d$-finiteness. Here, a proper connective dg algebra $A$ with $\gl A\leq d$ is said to be $\nu_d$-finite if its $d$-cluster category $\C_d(A)=(\per A/\nu_d)_{\triangle}$ is Hom-finite. We prove the following stronger statement, which is also of independent interest.

\begin{Thm}[Corollary \ref{nocyclepath}]
Let $A$ be a proper connective $\nu_d$-finite algebra with $\gl A\le d$. Then the degree-$d$ extension relation among the simple $A$-modules has no oriented cycle. Consequently, if $A=kQ$ is as above and is $\nu_d$-finite, then $Q$ has no cycle consisting entirely of arrows of degree $-d+1$.
\end{Thm}

It follows that the dg-quiver criterion applies at every vertex of a $\nu_d$-finite algebra. We also show that, under the silting--cluster-tilting correspondence of \cite{HI}, a silting mutation that remains $d$-silting induces the cluster tilting mutation of \cite{IYos}. This gives a conceptual link between preservation of global dimension and mutation in the associated cluster category.

We apply these results to higher hereditary algebras. For $d\geq1$, the notion of $d$-hereditary algebras is a generalization of path algebras to the case of global dimension is $d$ in the view point of higher Auslander-Reiten theory \cite{HIO}. They are considered as the most basic algebras among algebras of global dimension $d$ and possess beautiful properties generalizing those of path algebras. They consist of $d$-representation finite algebras and $d$-representation infinite algebras, which generalizes Dynkin/non-Dynkin dichotomy according to Gabriel's theorem \cite{HIO,IO11}. The two cases require different arguments, and we obtain a mutation theorem in each.

\begin{Thm}[Theorem \ref{mutdRI}]\label{mutdRIintro}
Let $A$ be a $d$-representation infinite algebra. Take $P\in\proj A$ and put $S:=\top P$. If $\pd_AS<d$ holds, then $\mu^-_P(A)$ is tilting and $\End_A(\mu^-_P(A))$ is a $d$-representation infinite algebra.
\end{Thm}

In the representation-finite case, the mutation must also respect the twisting automorphism in the fractionally Calabi--Yau structure. For a basic $l$-homogeneous $d$-representation finite algebra $A$, let $\phi$ be the automorphism from \cite[Proposition~4.3(a)]{HI11a} and put $\Phi:=-\otimes_A^{\mathbb L}A_\phi$, so that $\nu_d^{-l}\cong[d]\circ\Phi$.

\begin{Thm}[Theorem \ref{mutdRF}]
Write $A=P\oplus Q$. If
\[
 \Ext_A^d(\top P,Q)=0
 \qquad\text{and}\qquad
 \Phi(P)\cong P,
\]
then $\mu_P^-(A)$ is a tilting object and $\End_A(\mu_P^-(A))$ is again an $l$-homogeneous $d$-representation finite algebra.
\end{Thm}

The second condition ensures that the mutated object retains the same fractional Calabi--Yau symmetry; the first then allows the mutation criterion to force global dimension at most $d$. This is the point at which the representation-finite case differs essentially from the representation-infinite case.

These preservation theorems answer a question posed by Herschend--Iyama--Oppermann.

\begin{Ques}\cite[5.9]{HIO}\label{acyclicintro}
Is the Gabriel quiver of every higher hereditary algebra acyclic?
\end{Ques}

We give negative answers in both classes of the higher hereditary dichotomy.

\begin{Thm}[Examples \ref{counterHIOinfin}, \ref{counterHIOfin}]
We have the following examples of $2$-hereditary algebras with cycles.
\begin{enumerate}
\item A $2$-representation tame algebra with a $2$-cycle.
\item A $3$-homogeneous $2$-representation finite algebra with a cycle.
\end{enumerate}
\end{Thm}

\subsection*{Conventions}
Throughout this paper, $k$ denotes an arbitrary field. All algebras and categories are defined over $k$. For a dg algebra $A$, let $\D(A)$ denotes the unbounded derived category of right dg $A$-modules.

\subsection*{Use of AI}
The author used ChatGPT 5.6 Sol Pro solely in connection with Example \ref{counterHIOfin}. Specifically, when asked to search for a counterexample using Theorem \ref{mutdRF} as a guiding strategy, the model produced the example presented as Example \ref{counterHIOfin}. The author subsequently verified the example independently.

The author also used ChatGPT 5.6 Sol Extra High for language polishing and to improve the exposition of preliminary drafts. All AI-generated or AI-revised text was carefully reviewed and, where necessary, further revised by the author before being incorporated into the manuscript.

No AI tools were used to generate or modify any other mathematical arguments or results in this work.

\section*{Acknowledgements}
The author expresses his gratitude to Martin Herschend for suggesting a crucial idea and for many helpful discussions, especially concerning Subsection \ref{subsecmutdRF}. The author is also grateful to Norihiro Hanihara, Osamu Iyama and Nao Mochizuki for fruitful discussions. This work was supported by the WINGS-FMSP program at the Graduate School of Mathematical Sciences, the University of Tokyo, and JSPS KAKENHI Grant Number JP25KJ0818.

\section{Preliminaries on silting mutation}

Let $\T$ be a Hom-finite Krull-Schmidt triangulated category.

\begin{Def}
An object $M\in\T$ is called
\begin{enumerate}
\item {\it presilting} if $\T(M,M[>\!\!0])=0$ holds.
\item {\it pretilting} if $\T(M,M[\neq0])=0$ holds.
\item {\it silting} if it is presilting and $\T=\thick M$ holds.
\item {\it tilting} if it is pretilting and $\T=\thick M$ holds.
\end{enumerate}
We write $\silt\T$ for the isomorphism class of silting objects of $\T$.
\end{Def}

This set $\silt\T$ has several rich structures. First, we can equip a partial order with $\silt\T$ as follows.

\begin{Def}\cite[2.10,2.11]{AI}
For $M,N\in\silt\T$, we define
\[M\geq N:\Leftrightarrow\T(M,N[>\!\!0])=0.\]
Then this $\geq$ defines a partial order on $\silt\T$.
\end{Def}

Next, we can do an operation called {\it silting mutation} to elements in $\silt\T$.

\begin{Def}\cite[2.30,2.31,2.34]{AI}
Take $M\in\silt\T$. Decompose $M=X\oplus Y$ so that $(\add X)\cap(\add Y)=0$. Take a left $(\add Y)$-approximation $M\to Y_0$ and extend it to an exact triangle $M\to Y_0\to N\dashrightarrow$ which is called an exchange triangle. Then we call $\mu_X^-(M):=N$ a left mutation of $M$. Then $\mu_X^-(M)\in\silt\T$ holds. Dually, we define a right mutation $\mu_X^+(M)$.
\end{Def}

In \cite[2.36]{AI}, it is proved that for any $M,N\in\silt\T$ with $M>N$, there exists an indecomposable direct summand $X$ of $M$
 such that $\mu_X^-(M)\geq N$. In their proof, they give a sufficient condition for a direct summand $X$ of $M$ to be $\mu_X^-(M)\geq N$ (see also \cite[2.6]{HI}). In the following theorem, we show that this condition is also necessary.
 
\begin{Thm}\label{mutord}
Take $M,N\in\silt\T$ with $M\geq N$. Take a direct sum $M=X\oplus Y$ with $\add X\cap\add Y=0$. Let $M_0\xrightarrow{p}N$ be a minimal right $(\add M)$-approximation. Then the following conditions are equivalent.
\begin{enumerate}
\item $\mu^-_X(M)\geq N$
\item $M_0\in\add Y$
\end{enumerate}
\end{Thm}
\begin{proof}
Take a left $(\add Y)$-approximation $d\colon X\to Y_0$ and extend it to a triangle $X\xrightarrow{d}Y_0\to X'\dashrightarrow$. Then we have $\mu^-_X(M)=X'\oplus Y$.

(2)$\Rightarrow$(1) We give a complete proof for the convenience of the reader though it is parallel to that of \cite[2.36]{AI}. We have only to prove $\T(X',N[1])=0$. We have a long exact sequence
\[\T(Y_0,N)\to\T(X,N)\to\T(X',N[1])\to\T(Y_0,N[1]).\]
Since $\T(Y_0,N[1])=0$, it is enough to show that $\T(Y_0,N)\to\T(X,N)$ is surjective. Take $a\colon X\to N$. Then there exists $b\colon X\to M_0$ with $a=(X\xrightarrow{b}M_0\xrightarrow{p}N)$. Since $M_0\in\add Y$, $b$ factors through $d\colon X\to Y_0$.

(1)$\Rightarrow$(2) We prove every morphism $f\colon X\to M_0$ is a radical morphism. Since $\T(X'[-1],N)=0$, there exists $g\colon Y_0\to N$ with $pf=gd$. Then there exists $h\colon Y_0\to M_0$ with $g=ph$. Extend $p$ to a triangle $N'\xrightarrow{i}M_0\xrightarrow{p}N\dashrightarrow$. Since $p(f-hd)=0$, there exists $e\colon X\to N'$ with $f-hd=ie$. Since $d$ and $i$ are radical morphisms, we win.
\end{proof}

This characterization will play an essential role in our proof of Theorem \ref{dsiltmut}, \ref{mutdsilt}.

\section{$d$-silting objects}

Assume a triangulated category $\T$ satisfies the following conditions.

\begin{enumerate}
\item[(T0)] $\T$ is Hom-finite and Krull-Schmidt.
\item[(T1)] $\T$ has a Serre functor $\nu\curvearrowright\T$.
\end{enumerate}

We see that these conditions leads to a certain finiteness condition which corresponds to the properness of dg algebras.

\begin{Lem}\label{triproper}
Let $\T$ be a triangulated category satisfying (T0) and (T1). Then for any $X,Y\in\T$, we have $\T(X,Y[n])=0$ for $|n|\gg0$.
\end{Lem}
\begin{proof}
Take $X,Y\in\T$. By the existence of a silting object, we have $\T(X,Y[n])=0$ for $n\gg0$ \cite[2.4]{AI}. By the Serre duality, for $n\gg0$, we have
\[\T(X,Y[-n])\cong D\T(Y,\nu X[n])=0.\qedhere\]
\end{proof}

We put $\nu_d:=\nu\circ[-d]\curvearrowright\T$. We recall the definition of $d$-silting objects introduced by \cite{HI}.

\begin{Def}\cite[1.1]{HI}
Let $d\in\mathbb{Z}$ and $M\in\silt\T$. $M$ is called {\it $d$-silting} if $M\geq\nu_d^{-1}M$ holds. A tilting object which is $d$-silting is called {\it $d$-tilting}. Write
\[\silt^d\T:=\{M\in\silt\T\colon d\text{-silting}\}.\]
\end{Def}

Observe that for any silting object $M\in\T$, there must exist $d\in\mathbb{Z}$ such that $M$ is $d$-silting.

\begin{Ex}
\begin{enumerate}
\item If $A$ is a finite dimensional Iwanaga-Gorenstein $k$-algebra, then $A\in\per A$ is $d$-silting if and only if $\id_A A\leq d$ holds.
\item If $A$ is a $d$-selfinjective dg $k$-algebra in the sense of \cite{Jin}, then $A\in\per A$ is $(-d+1)$-silting.
\end{enumerate}
\end{Ex}

The following proposition states that the silting mutation at a direct summand only raises the dimension by at most one. In fact, this is a special case of Theorem \ref{mutdsilt}.

\begin{Prop}
Let $M\in\T$ be a $d$-silting object. Take a direct summand $N$ of $M$. Then $\mu^-_N(M)\in\T$ is a $(d+1)$-silting object.
\end{Prop}
\begin{proof}
Take an exchange triangle $M\to M'\to\mu^-_N(M)\dashrightarrow$. Then we have $\T(M,\nu^{-1}\mu^-_N(M)[>d])=0$. Thus we obtain $\T(\mu^-_N(M),\nu^{-1}\mu^-_N(M)[>d+1])=0$.
\end{proof}

We furthermore assume the following condition for $M\in\silt\T$. This condition correspondences to the finiteness of the global dimension.

\begin{enumerate}
\item[(T2)] $M$ admits a right adjacent t-structure $(\T_M^{\leq0}:=M[<0]^\perp,\T_M^{\geq0}:=M[>\!\!0]^\perp)$.
\end{enumerate}

Let $\H_M:=\T_M^{\leq0}\cap\T_M^{\geq0}$ be the heart of this t-structure. Observe that by combining with Lemma \ref{triproper}, we can say the following.

\begin{Lem}\label{fincoh}
Let $\T$ be a triangulated category satisfying (T0), (T1) and (T2). Then for any $X\in\T$, there exist integers $m<n$ such that $X\in\T_M^{\geq m}\cap\T_M^{\leq n}=\H_M[-n]*\H_M[-n+1]*\cdots\H_M[-m]$ holds.
\end{Lem}

Under these preparations, we can define the projective dimension of objects in $\T_{M}^{\leq0}$.

\begin{Propdef}\label{projdim}
For $T\in\T_M^{\leq0}$ and $d\geq0$, the following conditions are equivalent.
\begin{enumerate}
\item $T\in\add M*\add M[1]*\cdots*\add M[d]$
\item For any $H\in\H_M$, we have $\T(T,H[>d])=0$.
\end{enumerate}
If these conditions are satisfied, we write $\pd_MT\leq d$.
\end{Propdef}
\begin{proof}
(1)$\Rightarrow$(2) is obvious. We prove (2)$\Rightarrow$(1). Observe that by Lemma \ref{fincoh}, (2) is equivalent to that $\T(T,U)=0$ holds for all $U\in\T_M^{<-d}$. Consider a right $(\add M)$-approximation $M_0\to T$ and extend it to a triangle $T'\to M_0\to T\dashrightarrow$. Then by the long exact sequence induced by applying $\T(M,-)$ to this triangle, we obtain $T'\in\T_M^{\leq0}$. 

First, consider the case of $d=0$. Then since $T'[1]\in\T_M^{<0}$, we have $\T(T,T'[1])=0$. Thus we obtain $T\in\add M$. Next, consider the case of $d>0$. Then we can check $\T(T',H[>d-1])=0$ holds for any $H\in\H_M$. Then by inductive arguments, we have $T'\in\add M*\add M[1]*\cdots*\add M[d-1]$. Thus we obtain $T\in\add M*\add M[1]*\cdots*\add M[d]$.
\end{proof}

Finally, we see the following proposition which characterizes the $d$-siltingness.

\begin{Prop}\label{gldimd}
The following conditions are equivalent.
\begin{enumerate}
\item $M\in\silt^d\T$
\item $\pd_MH\leq d$ holds for all $H\in\H_M$.
\item $\T(\H_M,\H_M[>d])=0$
\item $\H_M\subseteq(\add M)*(\add M[1])*\cdots*(\add M[d])$
\item $\nu_d(\T_M^{\geq0})\subseteq\T_M^{\geq0}$
\item $\nu_d^{-1}(\T_M^{\leq0})\subseteq\T_M^{\leq0}$
\end{enumerate}
\end{Prop}
\begin{proof}
(2)$\Leftrightarrow$(3)$\Leftrightarrow$(4) follows from Proposition \ref{projdim}. We see (1)$\Leftrightarrow$(5)$\Leftrightarrow$(6). Observe that (1) is equivalent to $\nu_d^{-1}M\in\T_M^{\leq0}$. Since $\T_M^{\leq0}=\bigcup_{l\geq0}(\add M*\add M[1]*\cdots*\add M[l])$, this implies (1)$\Leftrightarrow$(6). (5)$\Leftrightarrow$(6) follows from Lemma \ref{fincoh}. See also \cite[4.1]{HI}. Next, we show (3)$\Rightarrow$(1). Remark that by Lemma \ref{fincoh}, (3) is equivalent to that $\T(X,Y)=0$ holds for every $X\in\T_M^{\geq0}$ and $Y\in\T_M^{<-d}$. By the Serre duality, we can check $\nu M\in\T_M^{\geq0}$. Thus we obtain $\T(\nu M,M[>d])=0$. 

Finally, we show (6)$\Rightarrow$(3). Take $H,H'\in\H_M$. Observe that we have $\T(H,\nu M[>\!\!0])\cong D\T(M[>\!\!0],H)=0$. Since $\nu M\in\silt\T$, there exists some $n\geq0$ such that $H\in\add\nu M[-n]*\cdots*\add\nu M[-1]*\add\nu M$ holds. Here, since $\nu_d^{-1}H'\in\T_M^{\leq0}$ holds by (6), we have $\T(\nu M,H'[>d])\cong\T(M,\nu_d^{-1}H'[>\!\!0])=0$. Therefore we obtain $\T(H,H'[>d])=0$.
\end{proof}

\section{Dg algebras}

\subsection{Global dimension}

In this subsection, we introduce the global dimension of locally finite connective dg algebras. First, we introduce basic terminologies.

\begin{Def}
A dg $k$-algebra $A$ is called
\begin{enumerate}
\item {\it locally finite} if $\dim_kH^nA<\infty$ holds for each $n\in\mathbb{Z}$.
\item {\it proper} if $\sum_{n\in\mathbb{Z}}\dim_kH^nA<\infty$ holds.
\item {\it connective} if $H^{>0}A=0$ holds.
\end{enumerate}
\end{Def}

Remark that $A$ is connective if and only if $A\in\per A$ is a silting object. Next, we introduce several subcategories of the derived category $\D(A)$ of a locally finite connective dg algebra $A$.

\begin{Def}
Let $A$ be a locally finite connective dg algebra.
\begin{enumerate}
\item $\per A:=\thick A\subseteq\D(A)$
\item $\pvd A:=\{M\in\D(A)\mid\sum_{n\in\mathbb{Z}}\dim_kH^nM<\infty\}$
\item $\D_{fd}(A):=\{M\in\D(A)\mid\dim_kH^nM<\infty\text{ holds for each }n\in\mathbb{Z}\}$
\item $\D^{\leq0}(A):=\{M\in\D(A)\mid H^{>0}M=0\}$
\item $\D^{\geq0}(A):=\{M\in\D(A)\mid H^{<0}M=0\}$
\end{enumerate}
\end{Def}

We also write $\D^{\leq0}_{fd}(A):=\D_{fd}(A)\cap\D^{\leq0}(A)$. Now, for a locally finite connective dg algebra $A$, we define the projective dimension of objects in $\D^{\leq0}_{fd}(A)$ and the global dimension of $A$.

\begin{Propdef}\label{gldimdg}
Let $A$ be a locally finite connective dg algebra. For $T\in\D^{\leq0}_{fd}(A)$ and $d\geq0$, the following conditions are equivalent.
\begin{enumerate}
\item $T\in\add A*\add A[1]*\cdots*\add A[d]$
\item For any $H\in\mod H^0A\subseteq\D(A)$, we have $\D(A)(T,H[>d])=0$.
\end{enumerate}
If these conditions are satisfied, then we write $\pd_AT\leq d$. If $\pd_AH\leq d$ holds for every $H\in\mod H^0A$, then we write $\gl A\leq d$. If such $d$ does not exist, then we write $\gl A=\infty$.
\end{Propdef}
\begin{proof}
This can be shown in the same way as Proposition-Definition \ref{projdim}. One thing where we have to be careful is to prove $\D(A)(T,X)=0$ for every $X\in\D^{<0}_{fd}(A)$ under the condition that $T$ satisfies $\D(A)(T,H[>\!\!0])=0$ for every $H\in\mod H^0A$. Since the system
\[H^n(\tau^{\ge-1}X)\longleftarrow H^n(\tau^{\ge-2}X)\longleftarrow\cdots\]
satisfies the Mittag-Leffler condition for each $n\in\mathbb{Z}$, we have
\[X\cong\holim(\tau^{\ge-1}X\longleftarrow\tau^{\ge-2}X\longleftarrow\cdots).\]
Since $\dim_k\D(A)(T[1],\tau^{\ge-n}X)<\infty$ holds for each $n\ge1$ by \cite[3.10]{Fus}, the system
\[\D(A)(T[1],\tau^{\ge-1}X)\longleftarrow \D(A)(T[1],\tau^{\ge-2}X)\longleftarrow\cdots\]
satisfies the Mittag-Leffler condition. Thus we have
\begin{align*}
\D(A)(T,X)&=\D(A)(T,\holim(\tau^{\ge-1}X\longleftarrow\tau^{\ge-2}X\longleftarrow\cdots)) \\
&\cong\lim(\D(A)(T,\tau^{\ge-1}X)\longleftarrow\D(A)(T,\tau^{\ge-2}X)\longleftarrow\cdots) \\
&=\lim(0\longleftarrow0\longleftarrow\cdots) \\
&=0.\qedhere
\end{align*}
\end{proof}

\begin{Rem}
In terms of Koszul duality \cite{Fus}, $\gl A<\infty$ holds if and only if its Koszul dual $A^!=\REnd_A(S_A)$ is proper where $S_A=\top H^0A\in\mod H^0A\subseteq\D(A)$. 
\end{Rem}

Then we can check the left-right symmetry of the global dimension.

\begin{Prop}
For a locally finite connective dg algebra $A$, we have
\[\gl A=\gl A^{\op}.\]
\end{Prop}
\begin{proof}
It is enough to show that $\gl A\le d$ implies $\gl A^{\op}\le d$. Take $H,H'\in\mod H^0A^{\op}\subseteq\pvd A^{\op}$. Since $D\colon\pvd A\to\pvd A^{\op}$ gives a duality, we have
\[\D(A^{\op})(H,H'[>d])\cong\D(A)(DH',DH[>d])=0.\qedhere\]
\end{proof}

We can characterize the finiteness of global dimension in the following way. Remark that $\pvd A\supseteq\per A$ holds if $A$ is proper.

\begin{Prop}\label{gldimfin}
For a locally finite connective dg algebra $A$, the following conditions are equivalent.
\begin{enumerate}
\item $\gl A<\infty$
\item $\pvd A\subseteq\per A$
\end{enumerate}
\end{Prop}
\begin{proof}
(1)$\Rightarrow$(2) Since $\pvd A=\thick(\mod H^0A)$, the assertion follows.

(2)$\Rightarrow$(1) Since $A$ is connective, for any $X,Y\in\per A$, there exists $d\in\mathbb{Z}$ such that $\D(A)(X,Y[>d])=0$. Since any object in $\mod H^0A$ can be written as a filtration of simple objects, whose number is finite, we can take $d\geq0$ such that for every $H,H'\in\mod H^0A$, we have $\D(A)(H,H'[>d])=0$.
\end{proof}

We see that if $A$ is proper, connective and $\gl A<\infty$, then $\per A$ admits a Serre functor.

\begin{Prop}\label{serreforproperdg}
Let $A$ be a proper connective dg algebra with $\gl A<\infty$. Then the auto-equivalence $\nu:=-\otimes_A^\mathbb{L}DA$ of $\per A$ is a Serre functor.
\end{Prop}
\begin{proof}
It is well-known that for $X,Y\in\per A$, we have $\RHom_A(X,Y)\cong D\RHom_A(Y,X\otimes_A^\mathbb{L}DA)$. Since $DA\in\pvd A=\per A$ holds by Proposition \ref{gldimfin}, our functor $\nu\colon\per A\to\per A$ is well-defined. By the Serre duality, this $\nu$ is fully-faithful. Here, observe that $D\colon\per A\to\per A^{\op}$ gives a duality. Thus we have $\per A=\thick DA$, which implies that $\nu$ is essentially surjective.
\end{proof}

From these preparations, we can check that for a proper connective dg algebra $A$ with $\gl A<\infty$, $M:=A\in\T:=\per A$ satisfy the conditions (T0), (T1) and (T2). Moreover, the definitions of the projective dimension of objects in $\T_A^{\leq0}$ in Proposition-Definition \ref{projdim} and \ref{gldimdg} coincide. In this setting, we write $\silt A:=\silt(\per A)$ and $\silt^dA:=\silt^d(\per A)$. By using Proposition \ref{gldimd}, we can see that $A\in\silt^dA$ if and only if $\gl A\leq d$.

\begin{Prop}\label{dgsiltinggldim}
For a proper connective dg algebra $A$ with $\gl A<\infty$, the following conditions are equivalent.
\begin{enumerate}
\item $A\in\silt^dA$
\item $\gl A\leq d$
\item $\nu_d(\per A\cap\D^{\geq0}(A))\subseteq\D^{\geq0}(A)$
\item $\nu_d^{-1}(\per A\cap\D^{\leq0}(A))\subseteq\D^{\leq0}(A)$
\end{enumerate}
\end{Prop}
\begin{proof}
This follows immediately from Proposition \ref{gldimd}.
\end{proof}

\subsection{Dg path algebras}

A {\it dg path algebra} is a dg algebra whose underlying graded algebra is a path algebra of a graded quiver. In this subsection, we give proofs to some folklores on dg path algebras. Let $A:=kQ$ be a dg path algebra where $\#Q_0<\infty$ and $Q_1^{>0}=\emptyset$. For $i\in Q_0$, let $S_i:=ke_i$ be the right simple $H^0A$-module corresponding to $i$, which we view as right dg $A$-module. We have a natural surjection $\pi\colon e_iA\to S_i$.

\begin{Prop}
$\Ker\pi$ is cofibrant as a right dg $A$-module. Thus $C:={\rm Cone}(\Ker\pi\to e_iA)$ gives a cofibrant resolution of $S_i$.
\end{Prop}
\begin{proof}
We can easily see $\Ker\pi=\bigoplus_{\alpha\in Q_1,t(\alpha)=i}\alpha A$ as a right $A$-module, but not as a right {\it dg} $A$-module! Let $F_n:=\bigoplus_{\alpha\in Q_1,t(\alpha)=i,|\alpha|\geq-n}\alpha A\subseteq\Ker\pi$ be a right sub dg $A$-module for $n\geq0$. Consider the following filtration of $\Ker\pi$.
\[0=:F_{-1}\subseteq F_0\subseteq F_1\subseteq\cdots\subseteq\Ker\pi\]
Then $\bigcup_{n\geq0}F_n=\Ker\pi$ holds. Moreover, we have short exact sequences
\[0\to F_{n-1}\to F_n\to\bigoplus_{\alpha\in Q_1,t(\alpha)=i,|\alpha|=-n}\alpha A\to0\ (n\geq0)\]
of right dg $A$-modules, where we view $\bigoplus_{\alpha\in Q_1,t(\alpha)=i,|\alpha|=-n}\alpha A$ as a direct sum of right dg $A$-modules $\alpha A\cong e_{s(\alpha)}A[-|\alpha|]$. Therefore we can conclude that $\Ker\pi$ is cofibrant.
\end{proof}

We prove that the extension groups between simple objects can be computed by counting the numbers of arrows. Observe that when $d=0$ and $H^{<0}A=0$, then this result is classical.

\begin{Thm}\label{extsimpath}
Assume $\#Q^{-d}_1<\infty$ for each $d\geq0$ and $d\alpha\in kQ_{\geq2}$ for each $\alpha\in Q_1$. Then for $i,j\in Q_0$ and $d\geq0$, we have
\[\dim_k\Ext^{d+1}_A(S_i,S_j)=\#\{\alpha\colon j\to i\mid|\alpha|=-d\}.\]
\end{Thm}
\begin{proof}
We have a short exact sequence
\[0\to\mathscr{H}om_A(\Ker\pi[1],S_j)\to\mathscr{H}om_A(C,S_j)\to\mathscr{H}om_A(e_iA,S_j)\to0.\]
Here $\mathscr{H}om_A(e_iA,S_j)=S_je_i=\delta_{ij}S_i$ holds. In addition, being induced by $\Ker\pi\to e_iA$, the map $H^0\mathscr{H}om_A(e_iA,S_j)\to H^1\mathscr{H}om_A(\Ker\pi[1],S_j)$ is $0$. Thus we have
\[H^d\mathscr{H}om_A(\Ker\pi,S_j)\cong H^{d+1}\mathscr{H}om_A(\Ker\pi[1],S_j)\xrightarrow{\cong}H^{d+1}\mathscr{H}om_A(C,S_j)=\Ext^{d+1}_A(S_i,S_j).\]

For $m\geq0$, we have short exact sequence
\[0\to\mathscr{H}om_A\Big(\bigoplus_{\alpha\in Q_1,t(\alpha)=i,|\alpha|=-m}\alpha A,S_j\Big)\to\mathscr{H}om_A(F_m,S_j)\to\mathscr{H}om_A(F_{m-1},S_j)\to0.\]
Here $\mathscr{H}om_A(\bigoplus_{\alpha\in Q_1,t(\alpha)=i,|\alpha|=-m}\alpha A,S_j)=\bigoplus_{\alpha\colon j\to i\in Q_1,|\alpha|=-m}S_j[-m]$ holds. For $\alpha\in Q_1$ with $t(\alpha)=i$ and $|\alpha|=-m$, we define a chain map $\alpha A[-1]\to F_{m-1}$ as $\alpha a\mapsto(d\alpha)a$. Then we can see $F_m={\rm Cone}(\bigoplus_{\alpha\in Q_1,t(\alpha)=i,|\alpha|=-m}\alpha A[-1]\to F_{m-1})$. Moreover, by our assumption $d\alpha\in kQ_{\geq2}$, the induced map $\mathscr{H}om_A(F_{m-1},S_j)\to\mathscr{H}om_A(\bigoplus_{\alpha\in Q_1,t(\alpha)=i,|\alpha|=-m}\alpha A,S_j)$ is $0$. Thus we have
\[H^l\mathscr{H}om_A(F_m,S_j)\xrightarrow{\cong}H^l\mathscr{H}om_A(F_{m-1},S_j)\ (l\neq m)\text{ and}\]
\[0\to\bigoplus_{\alpha\colon j\to i\in Q_1,|\alpha|=-m}k\to H^m\mathscr{H}om_A(F_m,S_j)\to H^m\mathscr{H}om_A(F_{m-1},S_j)\to0\colon\text{exact}.\]
Therefore for $m\geq d$, we have
\[H^d\mathscr{H}om_A(F_m,S_j)\xrightarrow{\cong}H^d\mathscr{H}om_A(F_d,S_j)\cong\bigoplus_{\alpha\colon j\to i\in Q_1,|\alpha|=-d}k.\]

We have $\mathscr{H}om_A(\Ker\pi,S_j)=\lim_{m\geq0}\mathscr{H}om_A(F_m,S_j)$. Observe that each term of $\mathscr{H}om_A(F_m,S_j)$ is finite dimensional. Therefore Mittag-Leffler conditions hold appropriately and we obtain
\[H^d\mathscr{H}om_A(\Ker\pi,S_j)=\lim_{m\geq0}H^d\mathscr{H}om_A(F_m,S_j)=\bigoplus_{\alpha\colon j\to i\in Q_1,|\alpha|=-d}k.\qedhere\]
\end{proof}

As a corollary, first, we can show that our dg path algebras are locally finite.

\begin{Cor}
Assume $\#Q^{-d}_1<\infty$ for each $d\geq0$ and $d\alpha\in kQ_{\geq2}$ for each $\alpha\in Q_1$. If $H^0A$ is finite dimensional, then $A$ is locally finite.
\end{Cor}
\begin{proof}
Since $\pvd A=\thick\{S_i\mid i\in Q_0\}$, Theorem \ref{extsimpath} implies that $\pvd A$ is Hom-finite. Thus by \cite[3.10]{Fus}, we get the assertion.
\end{proof}

Second, we can give an explicit formula of the global dimension of dg path algebras. When $A$ is proper, then this recovers \cite[8.2]{HI}.

\begin{Cor}\label{gldimpath}
Assume $\#Q^{-d}_1<\infty$ for each $d\geq0$ and $d\alpha\in kQ_{\geq2}$ for each $\alpha\in Q_1$. In addition, we assume $\dim_k H^0A<\infty$. Then for $d\geq0$, the following conditions are equivalent.
\begin{enumerate}
\item $\gl A\leq d$
\item $Q_1^{\leq-d}=\emptyset$
\end{enumerate}
In particular, $\gl A<\infty$ holds if and only if $Q$ is a finite quiver.
\end{Cor}

\section{Silting mutations preserving global dimension}

Assume a triangulated category $\T$ and $M\in\silt\T$ satisfy (T0), (T1), (T2) and the following condition.
\begin{enumerate}
\item[(T3)] $M\in\silt^d\T$
\end{enumerate}
Remark that we have a homological functor $H^0\colon\T\to\H_M$. In this section, we prove the following main theorems of this paper. The first one is characterizing when $\mu_X^{-}(M)\in\silt^d\T$ holds for $X\in\add M$.

\begin{Thm}\label{dsiltmut}
Assume a triangulated category $\T$ and $M\in\silt\T$ satisfy (T0), (T1), (T2) and (T3). Decompose $M=X\oplus X'$ with $(\add X)\cap(\add X')=0$ and put $S:=\top H^0X$. Then the following conditions are equivalent.
\begin{enumerate}
\item $\mu_X^{-}(M)\in\silt^d\T$
\item $\T(S,X'[d])=0$
\end{enumerate}
\end{Thm}

The second one is characterizing when $\mu_X^-(M)\geq\nu_d^{-1}M$ holds, which is a slightly stronger condition than $\mu_X^{-}(M)\in\silt^d\T$, for $X\in\add M$. As we will see, this characterization can be easily checked in terms of dg quivers.

\begin{Thm}\label{mutdsilt}
Assume a triangulated category $\T$ and $M\in\silt\T$ satisfy (T0), (T1), (T2) and (T3). Take $X\in\add M$ and put $S:=\top H^0X$. Then the following conditions are equivalent.
\begin{enumerate}
\item $\mu_X^-(M)\geq\nu_d^{-1}M$
\item $\pd_MS<d$
\item $\mu_X^-(M)\in\silt^d\T$ and $\T(S,S[d])=0$ holds.
\end{enumerate}
\end{Thm}

Towards these theorems, first, we exhibit a sequence of exact triangles which plays the same role as minimal injective resolutions.

\begin{Lem}\label{mininjres}
For $T=T_0\in\T_M^{\geq0}$, we have triangles
\[T_i\to\nu M^i\to T_{i+1}\dashrightarrow\ (i\geq0)\]
where $T_i\in\T_M^{\geq0}$ and the morphism $T_i\to\nu M^i$ is a minimal left $(\add\nu M)$-approximation.
\end{Lem}
\begin{proof}
We may assume $i=0$. Take a minimal left $(\add\nu M)$-approximation $T_0\to\nu M^0$ and extend it to a triangle $T_0\to\nu M^0\to T_1\dashrightarrow$. By applying $\T(-,\nu M)$ to this triangle, for $m>0$, we have an exact sequence
\[\T(\nu M^0,\nu M[m-1])\to\T(T_0,\nu M[m-1])\to\T(T_1,\nu M[m])\to\T(\nu M^0,\nu M[m]).\]
Observe that $\T(\nu M^0,\nu M[m])=0$ holds. If $m>1$, then since $\T(T_0,\nu M[m-1])\cong D\T(M[m-1],T_0)=0$, we have $\T(T_1,\nu M[m])=0$. If $m=1$, then since $\T(\nu M^0,\nu M)\to\T(T_0,\nu M)$ is surjective, we have $\T(T_1,\nu M[m])=0$. Therefore we obtain $\T(M[m],T_1)\cong D\T(T_1,\nu M[m])=0$. This means $T_1\in\T_M^{\geq0}$
\end{proof}

Next, we give an explicit formula of a minimal right $(\add M)$-approximation of $\nu_d^{-1}U$ where $U\in\T_M^{\leq0}$.

\begin{Lem}\label{constminrapp}
Let $U=U_0\in\T_M^{\leq0}\cap\T_M^{\geq-n}$ where $n\geq0$. Apply Lemma \ref{mininjres} to $T=U[-n]$ and obtain a sequence of exact triangles
\[U_i\to\nu M^i[n]\to U_{i+1}\xrightarrow{f_{i+1}}U_i[1]\ (i\geq0)\]
where $U_i\in\T_M^{\geq-n}$ and $U_i\to\nu M^i[n]$ is a minimal left $(\add\nu M[n])$-approximation.
\begin{enumerate}
\item For $0\leq i\leq n+d$, we have $\T(\T_M^{\geq-n},U_i[>n+d-i])=0$.
\end{enumerate}
Thus the triangle $U_{n+d}\to\nu M^{n+d}[n]\to U_{n+d+1}\dashrightarrow$ splits. This implies $U_{n+d}=\nu M^{n+d}[n]$.
\begin{enumerate}
\setcounter{enumi}{1}
\item The composition
\[f:=(\nu M^{n+d}[n]=U_{n+d}\xrightarrow{f_{n+d}}U_{n+d-1}[1]\xrightarrow{f_{n+d-1}[1]}\cdots\xrightarrow{f_1[n+d-1]}U_0[n+d])\]
gives a minimal right $(\add\nu M[n])$-approximation of $U_0[n+d]$. Therefore the composition
\[M^{n+d}=\nu^{-1}U_{n+d}[-n]\to \nu^{-1}U_{n+d-1}[-n+1]\to\cdots\to \nu^{-1}U_0[d]=\nu_d^{-1}U\]
gives a minimal right $(\add M)$-approximation of $\nu_d^{-1}U$.
\end{enumerate}
\end{Lem}
\begin{proof}
(1) By (1)$\Rightarrow$(3) of Proposition \ref{gldimd}, we have $\T(\T_M^{\geq-n},U[>n+d])=0$ since $U\in\T_M^{\leq0}$. Assume we have $\T(\T_M^{\geq-n},U_i[>n+d-i])=0$ for some $0\leq i<n+d$. By applying $\T(\T_M^{\geq-n},-)$ to the exact triangle $U_i\to\nu M^i[n]\to U_{i+1}\dashrightarrow$, for $m>n+d-i-1$, we obtain an exact sequence
\[\T(\T_M^{\geq-n},\nu M^i[n+m])\to\T(\T_M^{\geq-n},U_{i+1}[m])\to\T(\T_M^{\geq-n},U_i[m+1])).\]
By our assumption, we have $\T(\T_M^{\geq-n},U_i[m+1]))=0$. In addition, we have $\T(\T_M^{\geq-n},\nu M^i[n+m])\cong D\T(M^i[n+m],\T_M^{\geq-n})=0$ since $m>0$. Thus we obtain $\T(\T_M^{\geq-n},U_{i+1}[m])=0$.

(2) First, we show that $f\colon\nu M^{n+d}[n]=U_{n+d}\to U_0[n+d]$ is a right $(\add\nu M[n])$-approximation. By applying $\T(\nu M[n],-)$ to the exact triangle $U_i\to\nu M^i[n]\to U_{i+1}\dashrightarrow$ for $0\leq i<n+d$, we obtain an exact sequence
\[\T(\nu M[n],U_{i+1}[n+d-i-1])\to\T(\nu M[n],U_i[n+d-i])\to\T(\nu M[n],\nu M^i[2n+d-i])).\]
Since $\T(\nu M[n],\nu M^i[2n+d-i]))=0$, the map $\T(\nu M[n],U_{i+1}[n+d-i-1])\to\T(\nu M[n],U_i[n+d-i])$ is surjective. Thus the composition $\T(\nu M[n],U_{n+d})\to\T(\nu M[n],U_0[n+d])$ is surjective.

Second, we show that $f\colon U_{n+d}\to U_0[n+d]$ is right minimal. Take a morphism $g\colon U_{n+d}\to U_{n+d}$ such that $fg=f$ holds. Since $f_1[n+d-1]f_2[n+d-2]\cdots f_{n+d}(1_{U_{n+d}}-g)=0$, the morphism $f_2[n+d-2]\cdots f_{n+d}(1_{U_{n+d}}-g)\colon U_{n+d}\to U_1[n+d-1]$ factors through the morphism $\nu M^0[2n+d-1]\to U_1[n+d-1]$. Since $\T(U_{n+d},\nu M^0[2n+d-1])=0$, we obtain $f_2[n+d-2]\cdots f_{n+d}(1_{U_{n+d}}-g)=0$. By iterating this argument, we obtain $f_{n+d}(1_{U_{n+d}}-g)=0$. Thus the morphism $1_{U_{n+d}}-g$ factors through the morphism $\nu M^{n+d-1}[n]\to U_{n+d}$ which is a radical morphism. Thus $g$ is an isomorphism.
\end{proof}

Finally, we see how to compute the extension groups from simple objects in the heart.

\begin{Lem}\label{extfromsim}
Let $T=T_0\in\T_M^{\geq0}$. Apply Lemma \ref{mininjres} to $T$ and obtain a sequence of exact triangles
\[T_i\xrightarrow{a_i}\nu M^i\to T_{i+1}\dashrightarrow\ (i\geq0)\]
where $T_i\in\T_M^{\geq0}$ and $a_i$ is a minimal left $(\add\nu M)$-approximation. Then for $i\geq0$ and a simple object $S\in\H_M$, we have
\[\T(S,T[i])\cong\T(S,\nu M^i).\]
\end{Lem}
\begin{proof}
First, we show that the morphism $a_i\circ-\colon\T(S,T_i)\to\T(S,\nu M^i)$ is an isomorphism. Since $\T(S,T_{i+1}[-1])=0$, this map is injective. In what follows, we prove the surjectivity. Take a non-zero morphism $0\neq f\colon S\to\nu M_i$. Then $f$ factors through the morphism $H^0(\nu M_i)\to\nu M_i$. Thus we may view $S\subseteq H^0(\nu M_i)$ in the abelian category $\H_M$. Observe that we can also view $H^0(T_i)\subseteq H^0(\nu M_i)$. By octahedral axiom, there exist $E\in\T$ and the following commutative diagram of triangles.
\[\xymatrix{
S \ar[r] \ar@{=}[d] & H^0(\nu M_i) \ar[r] \ar[d] & H^0(\nu M_i)/S \ar@{.>}[r] \ar[d] & \\
S \ar[r] & \nu M^i \ar[r]^{b_i} \ar[d] & E \ar@{.>}[r] \ar[d] & \\
 & \tau_M^{>0}(\nu M_i) \ar@{=}[r] \ar@{.>}[d] & \tau_M^{>0}(\nu M_i) \ar@{.>}[d]\\
  & & 
}\]
By using octahedral axiom again, there exist $F\in\T$ and the following commutative diagram of triangles.
\[\xymatrix{
 & S \ar@{=}[r] \ar[d] & S \ar[d] \\
T_i \ar[r]^{a_i} \ar@{=}[d] & \nu M^i \ar[r] \ar[d]^{b_i} & T_{i+1} \ar@{.>}[r] \ar[d] & \\
T_i \ar[r]_{b_ia_i} & E \ar[r] \ar@{.>}[d] & F \ar@{.>}[r] \ar@{.>}[d] & \\
 & &
}\]
Then by the right most vertical triangle, we have an exact sequence
\[0\to H^{-1}(F)\to S\to H^0(T_{i+1})\]
in $\H_M$. Here, suppose that $S\cap H^0(T_i)=0$ holds as a subobject of $H^0(\nu M_i)$. Since we have an exact sequence $0\to H^0(T_i)\to H^0(\nu M^i)\to H^0(T_{i+1})$ in $\H_M$, this means that the morphism $S\to H^0(T_{i+1})$ is monic in $\H_M$. Thus we obtain $H^{-1}(F)=0$. Then since $\T(F[-1],\nu M^i)\cong D\T(M^i,F[-1])=0$, there exists $c_i\colon E\to\nu M^i$ such that $c_i(b_ia_i)=a_i$ holds. Since $a_i$ is left minimal, $c_ib_i$ is an isomorphism. Thus $b_i$ is a section. Since we have a triangle $S\to\nu M^i\xrightarrow{b_i}E\dashrightarrow$, this means that $S$ is a direct summand of $E[-1]$, but this contradicts to $E\in\T_M^{\geq0}$. Therefore $S\cap H^0(T_i)\neq0$ holds. Since $S$ is simple in $\H_M$, we obtain $S\subseteq H^0(T_i)$ as a subobject of $H^0(\nu M_i)$. This means that there exists $g\colon S\to T_i$ such that $a_ig=f$ holds.

By applying $\T(S,-)$ to the triangle $T_i\to\nu M^i\to T_{i+1}\dashrightarrow$, for $m>0$, we have an exact sequence
\[\T(S,T_i[m-1])\to\T(S,\nu M^i[m-1])\to\T(S,T_{i+1}[m-1])\to\T(S,T_i[m])\to\T(S,\nu M^i[m]).\]
Observe that $\T(S,\nu M^i[>\!\!0])\cong D\T(M^i[>\!\!0],S)=0$. Thus $\T(S,T_{i+1}[m-1])\to\T(S,T_i[m])$ is an isomorphism for $m>1$. If $m=1$, since $\T(S,T_i)\to\T(S,\nu M^i)$ is an isomorphism, so is $\T(S,T_{i+1})\to\T(S,T_i[1])$. Therefore we obtain
\[\T(S,T[i])=\T(S,T_0[i])\cong\T(S,T_1[i-1])\cong\cdots\cong\T(S,T_i)\cong\T(S,\nu M^i).\qedhere\]
\end{proof}

By combining these lemmas, we obtain the following corollary.

\begin{Cor}\label{criprojcov}
Let $U\in\T_M^{\leq0}$ and take a minimal right $(\add M)$-approximation $M_0\to\nu_d^{-1}U$. Take $X\in\add M$ and put $S:=\top H^0X$. Then the following conditions are equivalent.
\begin{enumerate}
\item $(\add M_0)\cap(\add X)=0$
\item $\T(S,U[d])=0$
\end{enumerate}
\end{Cor}
\begin{proof}
By Lemma \ref{fincoh}, we can take $n\geq0$ such that $U\in\T_M^{\leq0}\cap\T_M^{\geq-n}$ holds. Apply Lemma \ref{mininjres} to $T=U[-n]$ and obtain a sequence of exact triangles
\[U_i\to\nu M^i[n]\to U_{i+1}\xrightarrow{f_{i+1}}U_i[1]\ (i\geq0)\]
where $U_i\in\T_M^{\geq-n}$ and $U_i\to\nu M^i[n]$ is a minimal left $(\add\nu M[n])$-approximation. Then by Lemma \ref{constminrapp}, we have $M_0\cong M^{n+d}$. Thus (1) is equivalent to $\T(M^{n+d},S)=0$ since $H^0M\in\H_M$ is projective. On the other hand, by Lemma \ref{extfromsim}, we have
\[\T(S,U[d])\cong\T(S,T[n+d])\cong\T(S,\nu M^{n+d})\cong D\T(M^{n+d},S).\]
Thus the assertion follows.
\end{proof}

Under these preparations, we can prove our main theorems.

\begin{proof}[Proof of Theorem \ref{dsiltmut}]
Take a left $(\add X')$-approximation $X\to X'_0$ and extend it to an exact triangle $X\to X'_0\to Y\dashrightarrow$. Then we have $\mu_X^-(M)=Y\oplus X'$. By Theorem \ref{mutord} and Corollary \ref{criprojcov}, (1) is equivalent to $\T(S,(Y\oplus X')[d])=0$. By applying $\T(S,-)$ to the triangle $X\to X'_0\to Y\dashrightarrow$, we have an exact sequence
\[\T(S,X'_0[d])\to\T(S,Y[d])\to\T(S,X[d+1]).\]
Assume (2) holds. Then we have $\T(S,X'_0[d])=0$. In addition, since $X\in\T_M^{\leq0}$ and $\pd_MS\leq d$, we have $\T(S,X[d+1])=0$. Thus we obtain $\T(S,Y[d])=0$. This proves the assertion.
\end{proof}

\begin{proof}[Proof of Theorem \ref{mutdsilt}]
(1)$\Leftrightarrow$(2) By Theorem \ref{mutord} and Corollary \ref{criprojcov}, (1) is equivalent to $\T(S,M[d])=0$. This is equivalent to (2) by Proposition-Definition \ref{projdim} since $\pd_MS\leq d$.

(1)\&(2)$\Rightarrow(3)$ By (2), we have $\T(S,S[d])=0$. Since $M\geq\mu_X^{-}(M)$, we have $\nu_d^{-1}M\geq\nu_d^{-1}\mu_X^{-}(M)$. By combining this with $\mu_X^-(M)\geq\nu_d^{-1}M$, we obtain $\mu_X^-(M)\geq\nu_d^{-1}\mu_X^-(M)$.

(3)$\Rightarrow$(2) We may assume that we have a decomposition $M=X\oplus X'$ with $(\add X)\cap(\add X')=0$. Put $S':=\top H^0X'$. By octahedral axiom, there exist $E\in\T$ and the following commutative diagram of triangles.
\[\xymatrix{
\tau_M^{<0}X' \ar@{=}[r] \ar[d] & \tau_M^{<0}X' \ar[d] \\
E \ar[r] \ar[d] & X' \ar[r] \ar[d] & S' \ar@{.>}[r] \ar@{=}[d] & \\
\rad H^0X' \ar[r] \ar@{.>}[d] & H^0X' \ar[r] \ar@{.>}[d] & S' \ar@{.>}[r] & \\
  & 
}\]
Applying $\T(S,-)$ to the triangle $E\to X'\to S'\dashrightarrow$, we obtain an exact sequence
\[\T(S,X'[d])\to\T(S,S'[d])\to\T(S,E[d+1]).\]
By Theorem \ref{dsiltmut}, we have $\T(S,X'[d])=0$. Since $E\in\T_M^{\leq0}$ by the leftmost vertical triangle in the commutative diagram, we have $\T(S,E[d+1])=0$. Thus we obtain $\T(S,S'[d])=0$. Combining this with $\T(S,S[d])=0$, the assertion follows.
\end{proof}
As an immediate corollary, we obtain the following.

\begin{Cor}\label{easychar}
Assume a triangulated category $\T$ and $M\in\silt^d\T$ satisfy (T0), (T1) and (T2). Take $X\in\add M$ and put $S:=\top H^0X$. If $\T(S,S[d])=0$ holds, then the following conditions are equivalent.
\begin{enumerate}
\item $\mu_X^-(M)\in\silt^d\T$
\item $\pd_MS<d$
\end{enumerate}
\end{Cor}

In terms of dg quivers, we can rephrase our results in the following way.

\begin{Cor}\label{easycharpath}
Let $A=kQ$ be a proper dg path algebra such that $Q$ is a finite graded quiver with $Q_1^{>0}=Q_1^{\leq-d}=\emptyset$. We assume $d\alpha\in kQ_{\geq2}$ holds for each $\alpha\in Q_1$. For $i\in Q_0$, if there is no loop of degree $-d+1$ at $i$, then the following conditions are equivalent.
\begin{enumerate}
\item $\mu_{e_iA}^-(A)\in\silt^dA$
\item There is no arrow of degree $-d+1$ whose sink is $i$.
\end{enumerate}
\end{Cor}
\begin{proof}
This follows immediately from Theorem \ref{extsimpath} and \ref{mutdsilt}.
\end{proof}

We remark that this result for dg path algebras can be deduced from the explicit recipe in \cite{Opp}, but our proof is more conceptual.

\begin{Ex}
Let $A:=k[1\xrightarrow{\alpha}2\xrightarrow{\beta}3]/(\beta\alpha)$ be a path algebra with relation. Then $\gl A=2$ holds and $A$ is quasi-equivalent to the dg path algebra of the following dg quiver.
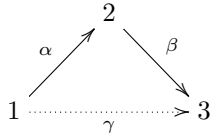
\[\xymatrix{
 & 2 \ar[dr]^\beta \\
1 \ar[ur]^\alpha \ar@{.>}[rr]_\gamma & & 3
}\]
Here, $\gamma$ denotes an arrow of degree $-1$ with $d\gamma=\beta\alpha$. Then by Corollary \ref{easycharpath}, $\mu_{e_iA}^-(A)\in\silt^2A$ holds if and only if $i=1,2$ since there is no loop of degree $-1$.
\end{Ex}


\section{Silting mutations for $\nu_d$-finite proper connective dg algebras}

In this section, we apply our main theorem to $\nu_d$-finite triangulated categories. First, we recall the definition of $\nu_d$-finiteness.

\begin{Def}\cite[4.7]{HI}
Let $\T$ be a triangulated category satisfying (T0) and (T1). We say that $\T$ is $\nu_d$-{\it finite} if for each $X,Y\in\T$, we have $\T(X,\nu_d^{-i}(Y)[\geq0])=0$ for $i\gg0$. A proper connective dg algebra $A$ with $\gl A<\infty$ is said to be $\nu_d$-{\it finite} if $\per A$ is $\nu_d$-finite.
\end{Def}

Observe that if $M\in\silt\T$ satisfies (T2), then $\T$ is $\nu_d$-finite if and only if for each $X\in\T$, we have $\nu_d^{\ll0}X\in\T_M^{\leq0}$.

\subsection{No cycles consisting of arrows of degree $-d+1$}

The following is our main theorem, which is of independent interest, in this subsection.

\begin{Thm}\label{nocycle}
Assume a triangulated category $\T$ and $M\in\silt\T$ satisfy (T0), (T1), (T2) and (T3). Moreover, we assume that $\T$ is $\nu_d$-finite. Then there exist no simple objects $S_1,\cdots,S_n, S_{n+1}=S_1\in\H_M$ such that $\T(S_i,S_{i+1}[d])\neq0$ holds for $1\leq i\leq n$. In particular, there exists no simple object $S\in\H_M$ such that $\T(S,S[d])\neq0$ holds.
\end{Thm}

To prove this theorem, we exhibit the following easy lemma.

\begin{Lem}\label{h0epi}
Assume a triangulated category $\T$ and $M\in\silt\T$ satisfy (T0), (T1) and (T2). Take an exact triangle $X\to Y\to Z\dashrightarrow$ with $Y,Z\in\T_M^{\leq0}$. Then the induced morphism $H^0(Y)\to H^0(Z)$ is epic in $\H_M$ if and only if $X\in\T_M^{\leq0}$ holds.
\end{Lem}

\begin{proof}[Proof of Theorem \ref{nocycle}]
Suppose that such simple objects $S_1,\cdots,S_n, S_{n+1}=S_1\in\H_M$ exist. By the Serre duality, we have a non-zero morphism $\nu_d^{-1}S_{i+1}\to S_i$ for $1\leq i\leq n$. Extend this to an exact triangle $X_i\to\nu_d^{-1}S_{i+1}\to S_i\dashrightarrow$. Observe that the induced morphism $H^0(\nu_d^{-1}S_{i+1})\to S_i$ is non-zero. Since $S_i\in\H_M$ is simple, this is epic. Thus by Lemma \ref{h0epi}, we have $X_i\in\T_M^{\leq0}$. By Proposition \ref{gldimd}, we have $\nu_d^{-m}X_i\in\T_M^{\leq0}$ for all $m\geq0$. Therefore again by Lemma \ref{h0epi}, the induced morphisms $H^0(\nu_d^{-m-1}S_{i+1})\to H^0(\nu_d^{-m}S_i)$ are all epic for $m\geq0$. This means that the compositions
\[\cdots\to H^0(\nu_d^{-2}S_{i+2})\to H^0(\nu_d^{-1}S_{i+1})\to S_i\]
are non-zero. Thus $H^0(\nu_d^{-m}S_i)\neq0$ holds for all $m\geq0$ and $1\leq i\leq n$. This contradicts to that $\T$ is $\nu_d$-finite.
\end{proof}

In terms of dg path algebras, we can rephrase our results in the following way.

\begin{Cor}\label{nocyclepath}
Let $A=kQ$ be a proper dg path algebra such that $Q$ is a finite graded quiver with $Q_1^{>0}=Q_1^{\leq-d}=\emptyset$. We assume $d\alpha\in kQ_{\geq2}$ holds for each $\alpha\in Q_1$. If $A$ is $\nu_d$-finite, then there exists no cycle consisting of arrows of degree $-d+1$. In particular, there is no loop of degree $-d+1$.
\end{Cor}

Thanks to Theorem \ref{nocycle}, we can restate Theorem \ref{mutdsilt} in the following simpler way.

\begin{Cor}
Assume a triangulated category $\T$ and $M\in\silt\T$ satisfy (T0), (T1), (T2) and (T3). Moreover, we assume that $\T$ is $\nu_d$-finite. Take an indecomposable direct summand $X$ of $M$ and put $S:=\top H^0X$. Then the following conditions are equivalent.
\begin{enumerate}
\item $\mu_X^-(M)\geq\nu_d^{-1}M$
\item $\pd_MS<d$
\item $\mu_X^-(M)\in\silt^d\T$
\end{enumerate}
\end{Cor}

\subsection{Compatibility with cluster tilting mutations}

Recall from \cite{IYos} that for a triangulated category $\T$ and $d\geq1$, a subcategory $\U\subseteq\T$ is called $d$-{\it rigid} if  $\T(\U,\U[i])=0$ holds for $0<i<d$. It is called $d$-{\it cluster tilting} if it is functorially finite, $d$-rigid and $\T=\U*\U[1]*\cdots*\U[d-1]$. We write $d$-$\ctilt\T:=\{\U\subseteq\T\colon d\text{-cluster tilting}\}$. In \cite{IYos}, mutations of cluster tilting subcategories are introduced.

\begin{Def}\cite[2.5,5.1]{IYos}
Let $\T$ be a triangulated category satisfying (T1) and (T2). For $\U\in d$-$\ctilt\T$ and a functorially finite subcategory $\D\subseteq\U$ with $\nu_d(\D)=\D$, define
\[\mu^-(\U;\D):=(\D*\U[1])\cap{}^\perp\D[1].\]
Then $\mu^{-}(\U;\D)\in d$-$\ctilt\T$ holds.
\end{Def}

On the other hand, in \cite{HI}, the following theorem, called silting-CT correspondence, is proved.

\begin{Thm}\cite[4.8]{HI}
Assume a triangulated category $\T$ and $M\in\silt\T$ satisfy (T0), (T1) and (T2). Moreover, we assume that $\T$ is $\nu_d$-finite. Then we have the following map.
\[\silt^d\T\to d\text{-}\ctilt\T;N\mapsto\U_d(N):=\add\{\nu_d^iN\mid i\in\mathbb{Z}\}\]
\end{Thm}

We prove the following compatibility between cluster tilting mutations and our silting mutations preserving global dimension. Compare this with \cite[4.25]{HI}.

\begin{Thm}
Assume a triangulated category $\T$ and $M\in\silt\T$ satisfy (T0), (T1), (T2) and (T3). Moreover, we assume that $\T$ is $\nu_d$-finite. Decompose $M=X\oplus X'$ so that $(\add X)\cap(\add X')=0$ holds. Put $\D:=\add\{\nu_d^iX'\mid i\in\mathbb{Z}\}\subseteq\U_d(M)$. If $\mu_X^-M\in\silt^d\T$ holds, then we have
\[\mu^-(\U_d(M);\D)=\U_d(\mu_X^-M).\]
\end{Thm}
\begin{proof}
Take a left $(\add X')$-approximation $X\to X'_0$. Then it is enough to show that this morphism is also a left $\D$-approximation. Observe that $\T(X,\nu_d^{>0}X')=0$ holds. Extend $X\to X'_0$ to an exact triangle $X\to X'_0\to Y\dashrightarrow$. Since $\mu_X^-M=Y\oplus X'\in\silt^d\T$, we have $\T(Y,\nu_d^{-m}X'[1])=0$ for $m>0$. Thus the induced morphism $\T(X'_0,\nu_d^{-m}X')\to\T(X,\nu_d^{-m}X')$ is surjective.
\end{proof}

\section{Silting mutations for higher hereditary algebras}

In this section, we apply Theorems \ref{dsiltmut} and \ref{mutdsilt} to obtain new higher hereditary algebras from the given ones by silting mutations. As an application, we give counterexamples to \cite[5.9]{HIO}.

First, we recall the definitions of $d$-representation finite, $d$-representation infinite and $d$-hereditary algebras introduced by \cite{IO11,HIO}. These algebras are introduced to give a class of finite dimensional algebras of global dimension $d$ sharing beautiful Auslander--Reiten-theoretic properties with hereditary algebras.

\begin{Def}\cite[2.2]{IO11}\cite[2.7,3.2]{HIO}
Let $A$ be a finite dimensional algebra. Let $d\geq1$ and assume $\gl A\le d$. $A$ is called
\begin{enumerate}
\item {\it $d$-representation finite} if $\mod A$ admits a $d$-cluster tilting object.
\item {\it $d$-representation infinite} if
\[\nu_d^{-n}A\in\mod A\subseteq\per A\]
holds for all $n\geq0$.
\item {\it $d$-hereditary} if for every $i\in\mathbb{Z}$, the cohomology groups of $\nu_d^iA$ is concentrated in degree of multiples of $d$.
\end{enumerate}
\end{Def}

In \cite{HIO}, they gave the following dichotomy.

\begin{Thm}\cite[3.4]{HIO}
Let $A$ be a ring-indecomposable finite dimensional algebra. Then $A$ is $d$-hereditary if and only if it is either $d$-representation finite or $d$-representation infinite.
\end{Thm}

In \cite{HIO}, the following question was exhibited.

\begin{Ques}\cite[5.9]{HIO}\label{acyclic}
The quivers of higher hereditary algebras are acyclic.
\end{Ques}

This question has long been open since it is posed. In what follows, we give counterexamples to this question for both $d$-representation finite and $d$-representation infinite cases as an application of our main theorem.

\subsection{Silting mutations for higher representation infinite algebras}

First, we prove that the silting mutation of $d$-representation infinite algebra satisfying the equivalent conditions in Theorem \ref{mutdsilt} is again $d$-representation infinite.

\begin{Thm}\label{mutdRI}
Let $A$ be a $d$-representation infinite algebra. Take $P\in\proj A$. If $M:=\mu^-_P(A)\geq\nu_d^{-1}A$ holds, then $M$ is tilting and $\End_A(M)$ is a $d$-representation infinite algebra.
\end{Thm}
\begin{proof}
By Theorem \ref{mutdsilt}, $M\in\silt^dA$ holds. For $n\geq0$, observe that we have
\[\D(A)(M[>\!\!0],\nu_d^{-n}M)\cong D\D(A)(\nu_d^{-n}M,\nu M[>\!\!0]).\]
Thus the assertion holds if and only if $\nu_d^{-n}M\geq\nu M$ holds for all $n\geq0$. By the same argument, we have $\nu_d^{-n-1}A\geq\nu A$. Therefore we obtain
\[\nu_d^{-n}M\geq\nu_d^{-n-1}A\geq\nu A\geq\nu M.\qedhere\]
\end{proof}

According to this theorem, it is natural to conjecture the following.

\begin{Conj}\label{dRIdsilt}
Let $A$ be a $d$-representation infinite algebra and $T\in\silt^dA$. Then $T$ is tilting and $\End_A(T)$ becomes $d$-representation infinite.
\end{Conj}

Observe that this conjecture is obviously true for $d=1$.

By using Theorem \ref{mutdRI}, we can give a counterexample to Question \ref{acyclic}.

\begin{Ex}\label{counterHIOinfin}
We view the polynomial ring $S:=k[x,y,z]$ as a $\mathbb{Z}$-graded $k$-algebra by $\deg x=\deg y=1$ and $\deg z=2$. Put $\qmod^\mathbb{Z}S:=\mod^\mathbb{Z}S/\fl^\mathbb{Z}S$ and write $\O\in\qmod^\mathbb{Z}S$ as the image of $S$ . Then $\E:=\bigoplus_{i=0}^3\O(i)\in\qmod^\mathbb{Z}S$ is a tilting object of $\D^b(\qmod^\mathbb{Z}S)$ and $A:=\End_{\qmod^\mathbb{Z}S}(\E)\cong\End_S^\mathbb{Z}(\bigoplus_{i=0}^3S(i))$ is a $2$-representation infinite algebra of type $\widetilde{A}$ (see \cite[4.2]{Tom25b}). The dg quiver description of $A$ is the following where the dotted arrows represent arrows of degree $-1$ whose differential give the commutative relations.
\[\xymatrix{
\O \ar@<0.25ex>[r]^x \ar@<-0.25ex>[r]_y \ar@<0.25ex>[dr]^z \ar@{.>}@<0.25ex>[d] \ar@{.>}@<-0.25ex>[dr] \ar@{.>}@<-0.25ex>[d] & \O(1) \ar@<0.25ex>[d]^x \ar@<-0.25ex>[d]_y \ar@<-0.25ex>[dl]_z \ar@{.>}@<0.25ex>[dl] \\
\O(3) & \O(2) \ar@<0.25ex>[l]^x \ar@<-0.25ex>[l]_y
}\]
We write $e_i\in A$ the corresponding idempotents for $0\leq i\leq 3$. Then by Theorem \ref{mutdsilt} and \ref{extsimpath}, $\mu^-_{e_iA}(A)\geq\nu_2^{-1}A$ holds if and only if $i=0,1$. Thus by Theorem \ref{mutdRI}, $B_i:=\End_A(\mu^-_{e_iA}(A))$ is $2$-representation infinite for $i=0,1$. Now we investigate the case of $i=0$. Put $\mathfrak{m}:=(x,y,z)\subseteq S$ and consider the graded Koszul complex of a regular sequence $x,y,z\in S$.
\[0\to S\to S(1)^{\oplus2}\oplus S(2)\to S(2)\oplus S(3)^{\oplus2}\to S(4)\to(S/\mathfrak{m})(4)\to0\]
This yields the following exact sequence in $\qmod^\mathbb{Z}S$.
\[0\to\O\xrightarrow{\phi}\O(1)^{\oplus2}\oplus\O(2)\to\O(2)\oplus\O(3)^{\oplus2}\to\O(4)\to0\]
Then we can easily see that $\phi$ is a left $(\add\bigoplus_{i=1}^3\O(i))$-approximation. Thus we have
\[B_0\cong\End_{\qmod^\mathbb{Z}S}\big(\Cok\phi\oplus\bigoplus_{i=1}^3\O(i)\big).\]
Here, by considering the degree $-2$ part of the graded Koszul complex, we can say there exists non-zero homomorphism $\O(2)\to\E$. By taking the dual, we can also conclude that there exists non-zero homomorphism $\E\to\O(2)$. Thus $B_0$ has a cycle.

In fact, by using the recipe in \cite{Opp}, we can calculate the dg quiver of $B_0$ and $B_1$.
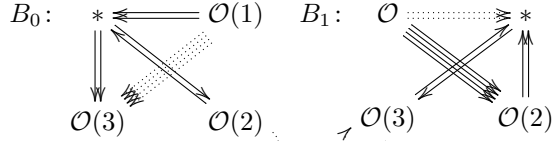
\[\begin{array}{c c}
B_0\colon\xymatrix{
{*} \ar@<0.25ex>[d] \ar@<-0.25ex>[d] \ar@<0.25ex>[dr] & \O(1) \ar@<0.25ex>[l] \ar@<-0.25ex>[l] \ar@<0.6ex>@{.>}[dl] \ar@<0.2ex>@{.>}[dl] \ar@<-0.2ex>@{.>}[dl] \ar@<-0.6ex>@{.>}[dl]\\
\O(3) & \O(2) \ar@<0.25ex>[ul] \ar@{.>}@(r,d)
}&B_1\colon\xymatrix{
\O \ar@<0.6ex>[dr] \ar@<0.2ex>[dr] \ar@<-0.2ex>[dr] \ar@<-0.6ex>[dr] \ar@{.>}@<0.25ex>[r] \ar@{.>}@<-0.25ex>[r] & {*} \ar@<0.25ex>[dl]\\
\O(3) \ar@<0.25ex>[ur] \ar@{.>}@(d,l) & \O(2) \ar@<0.25ex>[u] \ar@<-0.25ex>[u]
}\end{array}\]
Thus we can check that $B_0$ and $B_1$ have $2$-cycles directly. To understand them deeply, we draw the AR quiver of their $2$-preprojective components \cite[4.7]{HIO}.

\[
B_0\colon\xymatrix{
\O(1) \ar@<0.25ex>[dr] \ar@<-0.25ex>[dr] & \O(3) \ar@<0.6ex>[r] \ar@<0.2ex>[r] \ar@<-0.2ex>[r] \ar@<-0.6ex>[r] & \O(5) \ar@<0.25ex>[dr] \ar@<-0.25ex>[dr] & \O(7) \ar@<0.6ex>[r] \ar@<0.2ex>[r] \ar@<-0.2ex>[r] \ar@<-0.6ex>[r] & \O(9) \ar@<0.25ex>[dr] \ar@<-0.25ex>[dr] & \O(11) \ar@<0.6ex>[r] \ar@<0.2ex>[r] \ar@<-0.2ex>[r] \ar@<-0.6ex>[r] & \cdots\\
\O(2) \ar@<0.25ex>[r] \ar@/^-15pt/[rr] & {*} \ar@<0.25ex>[u] \ar@<-0.25ex>[u] \ar@<0.25ex>[l] & \O(6) \ar@<0.25ex>[r] \ar@/^-15pt/[rr] & {*}(4) \ar@<0.25ex>[u] \ar@<-0.25ex>[u] \ar@<0.25ex>[l] & \O(10) \ar@<0.25ex>[r] \ar@/^-15pt/[rr] & {*}(8) \ar@<0.25ex>[u] \ar@<-0.25ex>[u] \ar@<0.25ex>[l] & \cdots
}\]
\vspace{3mm}
\[B_1\colon\xymatrix{
\O \ar@<0.6ex>[r] \ar@<0.2ex>[r] \ar@<-0.2ex>[r] \ar@<-0.6ex>[r] & \O(2) \ar@<0.25ex>[d] \ar@<-0.25ex>[d] & \O(4) \ar@<0.6ex>[r] \ar@<0.2ex>[r] \ar@<-0.2ex>[r] \ar@<-0.6ex>[r] & \O(6) \ar@<0.25ex>[d] \ar@<-0.25ex>[d] & \O(8) \ar@<0.6ex>[r] \ar@<0.2ex>[r] \ar@<-0.2ex>[r] \ar@<-0.6ex>[r] & \O(10) \ar@<0.25ex>[d] \ar@<-0.25ex>[d] & \cdots \\
\O(3) \ar@<0.25ex>[r] \ar@/^-15pt/[rr] & {*} \ar@<0.25ex>[l] \ar@<0.25ex>[ur] \ar@<-0.25ex>[ur] & \O(7) \ar@<0.25ex>[r] \ar@/^-15pt/[rr] & {*}(4) \ar@<0.25ex>[l] \ar@<0.25ex>[ur] \ar@<-0.25ex>[ur] & \O(11) \ar@<0.25ex>[r] \ar@/^-15pt/[rr] & {*}(8) \ar@<0.25ex>[ur] \ar@<-0.25ex>[ur] \ar@<0.25ex>[l] & \cdots
}\]
\end{Ex}

\vspace{3mm}

\subsection{Silting mutations for higher representation finite algebras}\label{subsecmutdRF}

First, we recall the definition of $l$-homogeneous $d$-representation finite algebras introduced by \cite{HI11a}.

\begin{Def}\cite[1.2]{HI11a}
For $l\ge1$, a $d$-representation finite algebra $A$ is called {\it $l$-homogeneous} if
\[\nu_d^{-l+1}A\cong DA\]
holds in $\per A$.
\end{Def}

This notion can be rephrased by the property called twisted fractionally Calabi--Yau.

\begin{Thm}\cite[1.3]{HI11a}\label{homogtwfrCY}
For a finite-dimensional algebra $A$ and positive integers $d$ and $l$, the following conditions are equivalent.
\begin{enumerate}
\item $A$ is $l$-homogeneous $d$-representation finite.
\item $A$ is twisted $\frac{d(l-1)}{l}$-Calabi--Yau and $\gl A\le d$.
\end{enumerate}
\end{Thm}

Our starting point is a slight generalization of the implication (2)$\Rightarrow$(1) in Theorem \ref{homogtwfrCY}.

\begin{Prop}\label{twfrCYdRF}
Let $A$ be a proper connective dg algebra such that $\gl A\leq d$. Assume that there exists a triangle autoequivalence $\Phi\colon\per A\xrightarrow[\simeq]{}\per A$ such that $\nu^l\cong[d(l-1)]\circ\Phi$ for some $l\geq1$ and $\add\Phi(A)=\add A$. Then $H^{<0}A=0$. In particular, $A=H^0A$ is an $l$-homogeneous $d$-representation finite algebra.
\end{Prop}
\begin{proof}
Put $\T:=\per A$.  By Propositions \ref{serreforproperdg} and \ref{dgsiltinggldim}, $\T$ has the Serre functor $\nu=-\otimes_A^{\mathbb L}DA$, and $A\in\silt^d\T$.  Since every triangle autoequivalence commutes with a Serre functor up to a natural isomorphism, $\Phi\nu_d\cong\nu_d\Phi$.  Moreover, the assumed isomorphism of triangle functors gives
\[
 \nu_d^l=\nu^l[-dl]\cong\Phi[-d]
 \qquad\text{and hence}\qquad
 \nu_d^{-l}\cong\Phi^{-1}[d].
\]
Since $A$ is $d$-silting, we have $A\geq\nu_d^{-1}A$.  Applying powers of $\nu_d$ and using transitivity of the silting order, we obtain $\nu_d^qA\geq A$ for every $q\geq0$; in particular, $\nu_d^{l-1}A\geq A$.  For every $j>0$, Serre duality and the assumed isomorphism therefore yield
\begin{align*}
 H^{-j}A=\T(A,A[-j])
 &\cong D\T(A,\nu A[j]) \\
 &\cong D\T(\nu_d^{l-1}A,\nu_d^lA[d+j]) \\
 &\cong D\T(\nu_d^{l-1}A,\Phi(A)[j])=0.
\end{align*}
Indeed, the last equality follows from $\add\Phi(A)=\add A$ and $\nu_d^{l-1}A\geq A$.  Hence $H^{<0}A=0$.  By Theorem \ref{homogtwfrCY}, $A=H^0A$ is an $l$-homogeneous $d$-representation finite algebra.
\end{proof}

Combining with Theorem \ref{dsiltmut}, we obtain the following criterion for silting mutation to preserve $l$
-homogeneous $d$-representation finite algebras.

\begin{Thm}\label{mutdRF}
Let $A$ be a basic $l$-homogeneous $d$-representation finite algebra. By \cite[Proposition~4.3(a)]{HI11a}, choose a $k$-algebra automorphism $\phi$ of $A$ such that the triangle autoequivalence
\[
 \Phi:=-\otimes_A^{\mathbb L}A_\phi\colon\per A\xrightarrow{\sim}\per A
\]
satisfies $\nu_d^{-l}\cong[d]\circ\Phi$, where $A_\phi$ denotes the $A$--$A$-bimodule whose right action is twisted by $\phi$.  In particular, $\nu_d^{-l}A\cong A_\phi[d]$.  Take a decomposition $A=P\oplus Q$. Assume the following conditions are satisfied.
\begin{enumerate}
\item $\Ext^d_A(\top P,Q)=0$
\item $\Phi(P)\cong P$
\end{enumerate}
Then $\mu_P^-(A)\in\per A$ is a tilting object and $\End_A(\mu_P^-(A))$ is an $l$-homogeneous $d$-representation finite algebra.
\end{Thm}
\begin{proof}
Write $\T:=\per A$.  Let
\begin{equation}\label{mutdRFexchange}
 P\xrightarrow{f}Q_0\longrightarrow P'\longrightarrow P[1]
\end{equation}
be an exchange triangle in which $f$ is a minimal left $(\add Q)$-approximation, and put
\[
 M:=\mu_P^-(A)=P'\oplus Q.
\]
Taking $f$ to be minimal, $M$ is the basic representative of the mutated silting subcategory; see \cite[Theorem~2.31]{AI}.  By the defining property of $\Phi$, we have
\begin{equation}\label{mutdRFtwist}
 \nu_A^l\cong[d(l-1)]\circ\Phi^{-1},
\end{equation}
where $\nu_A$ denotes the Serre functor of $\T$.

We first show that $M$ is $\Phi$-stable.  As a right $A$-module, $A_\phi$ is isomorphic to $A$, and hence $\Phi(A)\cong A$.  The assumption $\Phi(P)\cong P$ and the Krull--Schmidt property therefore imply $\Phi(Q)\cong Q$.  Applying $\Phi$ to \eqref{mutdRFexchange}, we obtain a triangle
\[
 \Phi(P)\xrightarrow{\Phi(f)}\Phi(Q_0)\longrightarrow\Phi(P')
 \longrightarrow\Phi(P)[1].
\]
Here $\Phi(f)$ is again a minimal left $(\add Q)$-approximation: this follows from $\Phi(P)\cong P$, $\add\Phi(Q)=\add Q$, and the fact that an equivalence preserves minimal approximations.  Minimal left approximations are unique up to isomorphism, so the uniqueness of their cones gives $\Phi(P')\cong P'$.  Consequently, we obtain
\[
 \Phi(M)\cong M.
\]

The equality $\T(S,Q[d])=\Ext_A^d(S,Q)=0$ and Theorem~\ref{dsiltmut} show that
\[
 M\in\silt^d\T.
\]
Thus the assertions follow from Proposition \ref{twfrCYdRF}.
\end{proof}

\begin{Ex}\label{counterHIOfin}
Let $\Pi$ be the Jacobi algebra of a symmetric square-shaped quiver with potential of the follwoing form \cite[Subsection 9.3]{HI11b}.
\[
\begin{tikzcd}[row sep=1.8em, column sep=1.8em]
15 \arrow[dr]
  & 25 \arrow[l] \arrow[r]
  & 35 \arrow[d]
  & 45 \arrow[l] \arrow[d]
  & 55 \arrow[l] \\
14 \arrow[u] \arrow[d]
  & 24 \arrow[l] \arrow[d] \arrow[u]
  & 34 \arrow[l] \arrow[d] \arrow[ur]
  & 44 \arrow[l] \arrow[r]
  & 54 \arrow[u] \arrow[dl] \\
13 \arrow[d] \arrow[ur]
  & 23 \arrow[d] \arrow[l] \arrow[ur]
  & 33 \arrow[l] \arrow[r]
  & 43 \arrow[r] \arrow[u] \arrow[dl]
  & 53 \arrow[u] \arrow[dl] \\
12 \arrow[d] \arrow[ur]
  & 22 \arrow[l] \arrow[r]
  & 32 \arrow[r] \arrow[u] \arrow[dl]
  & 42 \arrow[r] \arrow[u] \arrow[d]
  & 52 \arrow[u] \arrow[d] \\
11 \arrow[r]
  & 21 \arrow[r] \arrow[u]
  & 31 \arrow[u]
  & 41 \arrow[l] \arrow[r]
  & 51 \arrow[ul]
\end{tikzcd}
\]
Then by \cite[9.9]{HI11b}, $\Pi$ is selfinjective, and its Nakayama automorphism is the half-turn
\[
 \sigma(i,j)=(6-i,6-j).
\]

Let $C$ be the cut indicated by the thick segments in the following diagram, and equip $\Pi$ with the induced $\mathbb{Z}_{\ge0}$-grading; see \cite[3.1]{HI11b}.
\[
\begin{tikzcd}[row sep=1.8em, column sep=1.8em]
15 \arrow[dr, line width=2pt, no head]
  & 25 \arrow[l] \arrow[r]
  & 35 \arrow[d]
  & 45 \arrow[l] \arrow[d]
  & 55 \arrow[l] \\
14 \arrow[u] \arrow[d]
  & 24 \arrow[l] \arrow[d] \arrow[u]
  & 34 \arrow[l, line width=2pt, no head] \arrow[d, line width=2pt, no head] \arrow[ur, line width=2pt, no head]
  & 44 \arrow[l] \arrow[r, line width=2pt, no head]
  & 54 \arrow[u] \arrow[dl] \\
13 \arrow[d, line width=2pt, no head] \arrow[ur, line width=2pt, no head]
  & 23 \arrow[d] \arrow[l] \arrow[ur]
  & 33 \arrow[l] \arrow[r]
  & 43 \arrow[r] \arrow[u] \arrow[dl]
  & 53 \arrow[u, line width=2pt, no head] \arrow[dl, line width=2pt, no head] \\
12 \arrow[d] \arrow[ur]
  & 22 \arrow[l, line width=2pt, no head] \arrow[r]
  & 32 \arrow[r, line width=2pt, no head] \arrow[u, line width=2pt, no head] \arrow[dl, line width=2pt, no head]
  & 42 \arrow[r] \arrow[u] \arrow[d]
  & 52 \arrow[u] \arrow[d] \\
11 \arrow[r]
  & 21 \arrow[r] \arrow[u]
  & 31 \arrow[u]
  & 41 \arrow[l] \arrow[r]
  & 51 \arrow[ul, line width=2pt, no head]
\end{tikzcd}
\]
The cut $C$ is invariant under $\sigma$. Then by \cite[2.3]{HI11a} and \cite[3.11]{HI11b}, $A:=\Pi_0$ is a $3$-homogeneous $2$-representation finite algebra. Put $P_{ij}:=e_{ij}A\in\proj A$ for $1\le i,j\le 5$. We perform the simultaneous left silting mutation at
\[P:=P_{43}\oplus P_{34}\oplus P_{23}\oplus P_{32}\]
simultaneously. An easy calculation, or Theorem \ref{extsimpath}, shows that
\[\pd_AS_{34}=\pd_AS_{32}=0\quad\text{and}\quad\pd_AS_{23}=\pd_AS_{43}=1\]
where $S_{ij}:=\top P_{ij}$. In addition, the half-turn exchanges $43$ with $23$ and $34$ with $32$. Thus the four vertices defining $P$ form a union of orbits of the Nakayama permutation, and the autoequivalence $\Phi$ in Theorem \ref{mutdRF} satisfies $\Phi(P)\cong P$.  Theorem \ref{mutdRF} therefore shows that
\[
 M:=\mu_P^-(A)
\]
is a tilting object in $\per A$ and that
\[
 B:=\End_A(M)
\]
is a $3$-homogeneous $2$-representation finite algebra.  We show that the Gabriel quiver of $B$ nevertheless contains an oriented cycle.

Write $A=P\oplus Q$. Reading off the minimal left $(\add Q)$-approximations from the quiver of $A$, we obtain the following objects of $\K^{[-1,0]}(\proj A)$:
\begin{align*}
 P'_{43}&:=\bigl[P_{43}\longrightarrow P_{33}\oplus P_{42}\oplus P_{54}\bigr],\\
 P'_{34}&:=\bigl[P_{34}\longrightarrow P_{12}\oplus P_{35}\oplus P_{44}\bigr],\\
 P'_{23}&:=\bigl[P_{23}\longrightarrow P_{12}\oplus P_{24}\oplus P_{33}\bigr],\\
 P'_{32}&:=\bigl[P_{32}\longrightarrow P_{22}\oplus P_{31}\oplus P_{54}\bigr].
\end{align*}
Here the left and right terms are placed in degrees $-1$ and $0$, respectively.  Hence
\[
 M=P'\oplus Q,
 \qquad
 P':=P'_{43}\oplus P'_{34}\oplus P'_{23}\oplus P'_{32}.
\]
To prove the existence of a cycle, it is enough to show that there are non-zero morphisms
\[P'_{32}\to P_{21}\to P_{11}\to P_{12}\to P'_{34}\to P_{45}\to P_{55}\to P_{54}\to P'_{32}.\]
Note that we claim that each morphism is non-zero, and do {\it not} claim that the composition of them is non-zero. By the symmetry with respect to $\sigma$, it suffices to construct non-zero morphisms
\[P'_{32}\to P_{21}\to P_{11}\to P_{12}\to P'_{34}.\]
Obviously, we have non-zero morphisms $P_{21}\to P_{11}\to P_{12}$.

We construct a non-zero morphism $P'_{32}\to P_{21}$.  Denote the four arrows around the lower-left square by
\[
 a\colon21\longrightarrow22,
 \qquad b\colon22\longrightarrow32,
 \qquad c\colon21\longrightarrow31,
 \qquad d\colon31\longrightarrow32.
\]
The cyclic derivative of the potential with respect to the cut arrow $32\to21$ yields a relation
\begin{equation}\label{counterHIOfiniterelation}
 ba-dc=0
\end{equation}
in $A$.  For a path $p\colon i\to j$, write $[p]\colon P_j\to P_i$ for the induced morphism of right projective $A$-modules.  The first two components of the differential of $P'_{32}$ are $[b]$ and $[d]$.  Therefore
\[
 u^0:=\bigl([a],-[c],0\bigr)
 \colon P_{22}\oplus P_{31}\oplus P_{54}\longrightarrow P_{21}
\]
satisfies $u^0d_{P'_{32}}=0$ by \eqref{counterHIOfiniterelation}, and hence defines a chain map $u\colon P'_{32}\longrightarrow P_{21}$. The map $u^0$ is nonzero.  Moreover, since $P_{21}$ is concentrated in degree $0$, no nonzero chain map of this form is null-homotopic.  Thus $u$ is nonzero in $\K^b(\proj A)$.

Finally, define a chain map $v\colon P_{12}\longrightarrow P'_{34}$ whose degree-zero component is the canonical inclusion
\[
 v^0=\begin{pmatrix}1_{P_{12}}&0&0\end{pmatrix}^{\mathsf T}
 \colon P_{12}\longrightarrow P_{12}\oplus P_{35}\oplus P_{44}.
\]
If $v$ were null-homotopic, there would exist a morphism $h\colon P_{12}\to P_{34}$ such that $v^0=d_{P'_{34}}h$.  Projecting this equality onto the $P_{12}$-summand would make the corresponding component $P_{34}\to P_{12}$ of $d_{P'_{34}}$ a split epimorphism, which is impossible.  Thus $v$ is nonzero in $\K^b(\proj A)$.
\end{Ex}

\begin{appendix}

\end{appendix}

\bibliographystyle{alpha} 
\bibliography{reference}

\end{document}